\documentclass[12pt, oneside]{scrartcl}
\usepackage[english]{babel}
\usepackage[T1]{fontenc}
\usepackage{amssymb}
\usepackage{amsmath}
\usepackage{hyperref}
\usepackage{amsthm}
\usepackage{cleveref}
\usepackage{tikz}
\newtheorem{theorem}{Theorem}[section]
\newtheorem{lemma}{Lemma}[section]
\usepackage{float}

\usepackage{tikz}
\usepackage{calc}
\usetikzlibrary{angles,quotes}
\usetikzlibrary{calc}
\usetikzlibrary{patterns}
\usetikzlibrary{graphs}
\usetikzlibrary{graphs.standard}
\usetikzlibrary{decorations.markings,shapes.geometric,positioning}
\title{A Tighter Bound on the 2-Distance Chromatic Number of Planar Graphs with Specific Maximum Degree}
\author{Sara Al Hajjar \footnote{Kalma Laboratory, 
Lebanese Univeristy, Beirut, Lebanon \\
Univ. Bordeaux, CNRS, Bordeaux INP, LaBRI, UMR 5800, F-33400, Talence,
France
}}
\date{}
\begin{document}

\maketitle
\begin{abstract}
 A $2$-distance $k$-coloring of a graph is a coloring of the vertices of the graph using $k$ colors such that any two vertices at distance two or less receive different colors. The 2-distance chromatic number of a graph $G$, denoted as $\chi_{2}(G)$, is the minimum integer $k$ such that $G$ has a 2-distance $k$-coloring. M. Krzyzinski \cite{7} proved that $\chi_2(G)\leq 3\Delta +4$ for planar graphs with maximum degree $\Delta \ge 6.$
  Later, N. Bousquet \cite{3} proved that $\chi_2(G) \leq 2\Delta +7$ for  planar graphs with $\Delta \ge 9$ and that  $\chi_2(G) \leq 21$ for planar graphs with  $\Delta \leq 6$. In this paper, we  prove that $\chi_2(G) \leq 3\Delta +2$ for  planar graphs  with a maximum degree $\Delta \ge 6$ hence improving the bound on $\chi_2(G)$ for $6 \leq \Delta \leq 8$.
\end{abstract}


\section{Introduction}
\begin{large}
Throughout this paper, we consider only finite simple graphs and use standard notations. The set of neighbors of a vertex $v$ in a graph $G$ is denoted by $N_G(v)$. The degree of a vertex $v$ in $G$ is the number of its neighbors and is denoted by $d_G(v)$. For brevity, we write $N(v)$ and $d(v)$ instead of $N_G(v)$ and $d_G(v)$, respectively. We denote by $\Delta(G)$ and $\delta(G)$ the maximum degree and minimum degree of a graph $G$, respectively. A vertex $v$ in $G$ is called a $k$-vertex, where $0 \leq k \leq \Delta(G)$, if $d(v) = k$. Moreover, a vertex $v$ in $G$ is called a $k^+$-vertex (resp., a $k^-$-vertex) if its degree is at least $k$ (resp., at most $k$). A $k$-neighbor (resp., $k^+$-neighbor, $k^-$-neighbor) of $v$ is a $k$-vertex (resp., $k^+$-vertex, $k^-$-vertex) adjacent to $v$.
The distance between two vertices $x$ and $y$ is the length of the shortest path connecting $x$ and $y$ in $G$ and is denoted by $d(x,y)$. For $i \ge 2$, the set $N_i(v)$ is defined as the set of all vertices of $G$ at distance at most $i$ from $v$ and we set $d_i(v) = |N_i(v)|$. For a set $S \subseteq V(G)$, we denote by $G[S]$ the subgraph of $G$ induced by the vertices in $S$.
A planar graph is a graph that can be drawn in the plane with no edge crossings. Such a drawing is called a plane graph or a planar embedding of the graph. When a planar graph is drawn without edge crossings, it partitions the plane into a set of regions, called faces. The set of faces of a planar graph $G$ is denoted by $F(G)$. Each face $f \in F(G)$ is bounded by a closed walk, called the boundary of $f$. The degree of a face $f$ is defined as the length of its boundary and is denoted by $d(f)$. A face $f$ in a planar graph $G$ is said to be incident with the vertices and edges on its boundary. Two faces of $G$ are said to be adjacent if their boundaries share a common edge. A face $f$ in a planar graph $G$ is called a $k$-face if $d(f)=k$. Moreover, a face $f$ in $G$ is called a $k^+$-face (resp., a $k^-$-face) if $d(f) \ge k$ (resp., $d(f) \le k$). \\

Let $G$ be a planar graph.  
A $(k,d)$-vertex in $G$ is a $k$-vertex incident to exactly $d$ $3$-faces.  
A $(k,d_1,d_2)$-vertex is a $k$-vertex incident to exactly $d_1$ $3$-faces and $d_2$ $4$-faces. 
A $(k,d^+)$-vertex is a $k$-vertex incident to at least $d$ $3$-faces. 
Similarly, a $(k,d^-)$-vertex is a $k$-vertex incident to at most $d$ $3$-faces. \\

A $2$-distance $k$-coloring of a graph $G$ is a mapping $\phi : V(G) \to \{1,2,\ldots,k\}$ such that $\phi(v_1) \neq \phi(v_2)$ whenever $d(v_1,v_2) \leq 2$, where $v_1$ and $v_2$ are any two vertices of $G$.  
The $2$-distance chromatic number of a graph $G$, denoted by $\chi_2(G)$, is the minimum integer $k$ such that $G$ admits a $2$-distance $k$-coloring.
\\ 

Several papers have studied Wegner's conjecture \cite{11} regarding the  $2$-distance chromatic number of planar graphs. Wegner conjectured the following:
\\ \\
\textbf{Wegner's Conjecture \cite{11}:} If $G$  is a planar graph with maximum degree $\Delta$, then 
$\chi_2(G) \leq 7$ if $\Delta=3$, $\chi_{2}$ ($G$ ) $\leq$ $\Delta + 5$ if $4 \leq \Delta \leq 7$ and $\chi_{2}(G)$ $\leq \lfloor\frac{3\Delta}{2}\rfloor+1$ if $\Delta \ge 8$. \\ \\
This conjecture remains largely unsolved. Thomassen \cite{9} proved
the conjecture for planar graphs with $\Delta=3$. 
For planar graphs with larger maximum degree, some upper bounds on the 2-distance chromatic number are already established. Agnarsson and Halldorsson \cite{1} showed that $\chi_2(G) \leq \lfloor\frac{ 9\Delta}{5}\rfloor +2 $ for planar graphs with maximum degree $ \Delta \ge 749$. Borodin et al. \cite{2} then improved the bound of $\chi_2(G)$ by proving that $\chi_2(G) \leq \lceil \frac{9\Delta}{5} \rceil+1$ for planar graphs with maximum degree $\Delta \ge 47$. Van de Heuvel and McGuinness \cite{10} showed that $\chi_2(G) \leq 2\Delta +25$ with no restriction on $\Delta$ while the bound $\chi_2(G) \leq \lceil \frac{5\Delta}{3}\rceil  + 78$ was proved by Molloy and Salavatipour \cite{8}. Zhu and Bu \cite{12} proved $\chi_2(G) \leq 5\Delta -7 $ when $\Delta \ge 6$ and $\chi_2(G) \leq 5\Delta - 9$ for $\Delta \ge 7 $ hence improving the bound of $\chi_2(G)$ for $6 \leq \Delta \leq 10$. Later, M. Krzyzinski \cite{7} proved that $\chi_2(G)\leq 3\Delta +4$ for $\Delta \ge 6.$ Finally, N. Bousquet \cite{3} proved that $\chi_2(G) \leq 2\Delta+7$ for $\Delta \ge 9$ improving the bound for $9 \leq \Delta \leq 31$ and $\chi_2(G) \leq 21$ for $\Delta \leq 6$ improving the bound for $\Delta =6$.
Moreover, Zhu and Bu \cite{12} showed that $\chi_2(G) \leq 20$ for planar graphs with maximum degree $\Delta \leq 5$. This bound was later reduced to 19 by Chen \cite{4} and to 18 by Aoki \cite{6}. Zou et al. \cite{13} then reduced the bound to 17 and finally Zakir Deniz \cite{5} reduced it to 16. \\

In this paper, we are going to prove that for a planar graph $G$ with maximum degree $\Delta \ge 6$, we have $\chi_2(G) \leq 3\Delta +2$ improving the bound of $\chi_2(G)$ proved in \cite{3} and \cite{7} for $6 \leq \Delta \leq 8$.

\section{Main Result:}

\begin{theorem} Let $G$  be a planar graph with maximum degree $\Delta \ge 6$, then $\chi_2(G)\leq 3\Delta +2$.
\end{theorem}

Our proof proceeds by contradiction. Assume that the statement is false and let $G$ be a minimal counterexample on $|E(G)|+|V(G)|$ not satisfying Theorem 2.1; $G$ is a planar graph with $\Delta(G) \ge 6$ but $\chi_2(G) > 3\Delta +2.$ We will prove that such a graph does not exist and hence Theorem 2.1 is true. \\

To handle such a graph $G$, we will use the following notion: 
We say that a graph $H$ is proper with respect to $G$ if $H$ is obtained from $G$ by deleting some vertices or edges and possibly adding edges, such that for every pair of vertices $v_1$ and $v_2$ in $V(G) \cap V(H)$ with $d_G(v_1,v_2) \leq 2$, we also have $d_H(v_1,v_2) \leq 2$, and $\Delta(H) \leq \Delta(G)$. Note that if $\phi$ is a $2$-distance coloring of such a graph $H$, then $\phi$ can be extended to the entire graph $G$, provided that each of the remaining uncolored vertices in $G$ has a safe color. \\
Our proof proceeds as follows: We present some structural results and forbidden configurations for the graph $G$. We then apply the discharging method together with Euler’s formula to derive a contradiction, showing that no counterexample to Theorem 2.1 exists. Hence, Theorem 2.1 holds. \\

By Euler's formula for planar graphs, we have the following equality:
\\
$$\sum_{v\in V(G)}(d(v)-4)+ \sum_{f \in F(G) }(d(f)-4) = -8$$
We assign an initial charge $d(v)-4$ to every vertex $v$ in $G$ and $d(f)-4$ to every face $f$ in $G$. We then assign appropriate discharging rules and  redistribute charges among vertices and faces, such that the final charge of each vertex
and each face is nonnegative. This would contradict the above equality since the sum of the charges should remain equal to $-8$. To ensure the existence of such a distribution, we must carefully analyze the structure of the graph $G$: In particular, the degrees and configurations of neighboring vertices and faces. This structural analysis is crucial to guarantee that the discharging process preserves non-negativity among all the vertices and faces of the graph.

\subsection{A Structural Analysis of Minimal Counterexample}
\begin{lemma} \label{l1} $G$ has no cut vertex. \end{lemma}
\begin{proof} Suppose that $G$ has a cut vertex $v$ and let $C_1, ..., C_t$ be the connected components of $G -v$, $t\ge 2$. Let $G_1 = C_1 \cup \{v\}$ and $G_2 = C_2\cup ... \cup C_t \cup \{v\}$. By definition of minimal counterexample, we have $\chi_2(G_i) \leq 3\Delta+2$ for $i=1,2$. Consider a $2$-distance ($3\Delta +2$)-coloring of $G_1$ and $G_2$ using the same colors for $G_1$ and $G_2$ such that $v$ receives the same color in $G_1$ and $G_2$. We first ensure that the neighbors of $v$ receive pairwise distinct colors. Indeed, such a coloring can be obtained by a color-switching argument, since $d(v) \leq \Delta$ and we have ($3\Delta + 2$) colors available. Then, by combining the colorings of $G_1$ and $G_2$, we obtain $(3\Delta+2)$-coloring of $G$, \mbox{a contradiction.} 
\end{proof}
We deduce that every face in $G$ is a cycle and every $k$-vertex is incident to exactly $k$ faces. 
\begin{lemma} \label{l2} Every vertex in $G$ is of degree at least 3.
\end{lemma}
\begin{proof} Suppose that there exists a vertex $v\in V(G)$ such that $d(v) \leq 2$. Let $N(v)=\{ v_1, v_2\}$. Let $G'$ be the graph obtained from $G$  after deleting $v$ and adding edge $v_1v_2$. $G'$ is proper with respect to $G$. By definition of minimal counterexample, we have $\chi_2(G') \leq 3\Delta +2$. Since $d_2(v) \leq 2\Delta < 3\Delta+2$, we color $v$ by a safe color in $G$ to obtain $\chi_2(G) \leq 3\Delta+2$, a contradiction. \end{proof} 
\begin{lemma} \label{l3} Let $v$ be a $3$-vertex. Then, we have the following:
\begin{enumerate}
    \item \label{l3p1} The neighbors of $v$ are of degree $\Delta$.
    \item \label{l3p2} $v$ is not incident to any $3$-face.
    \item \label{l3p3} $v$ is incident to at most one $4$-face.
\end{enumerate}
\end{lemma}
\begin{proof}  Let $v$ be a $3$-vertex and set $N(v)=\{v_1,v_2,v_3\}$.
\begin{enumerate}
\item Suppose that $v$ is adjacent to a $(\Delta-1)^-$-vertex.  W.l.o.g, suppose that $v_2$ is a $(\Delta-1)^-$-vertex. Let $G'$ be the graph obtained from $G$ after deleting $v$ and adding edges $v_1v_2$ and $v_2v_3$. $G'$ is proper with respect to $G$. By definition of minimal counterexample, we have $\chi_2(G') \leq 3\Delta +2$. Since $d_2(v) \leq 3\Delta -1 < 3\Delta+2$, we color $v$ by a safe color in $G$ to obtain $\chi_2(G) \leq 3\Delta +2$, a contradiction.  
\item Suppose that $v$ is incident to a $3$-face. W.l.o.g, suppose that the $3$-face incident to $v$ is $vv_1v_2.$ Let $G'$ be the graph obtained from $G$  after deleting $v$ and adding edge $v_1v_3$. $G'$ is proper with respect to $G$. By definition of minimal counterexample, we have $\chi_2(G') \leq 3\Delta +2.$ Since $d_2(v) < 3\Delta+2$, we color $v$ by a safe color in $G$ to obtain $\chi_2(G) \leq 3\Delta +2$, a contradiction. 
\item Suppose that $v$ is incident to two $4$-faces $f_1$ and $f_2$. W.l.o.g, suppose that $f_1 = vv_1xv_2$ and $f_2 = vv_2yv_3$ where $x \in N(v_1) \cap N(v_2)$ and $y \in N(v_2) \cap N(v_3)$. Let $G'$ be the graph obtained from $G$ after deleting $v$ and adding edge $v_1v_3$.
$G'$ is proper with respect to $G$. By definition of minimal counterexample, we have $\chi_2(G') \leq 3\Delta +2.$ Since $d_2(v)  < 3\Delta+2$, we color $v$ by a safe color in $G$ to obtain $\chi_2(G) \leq 3\Delta +2$, a contradiction.  
\end{enumerate} 
\end{proof}
\begin{lemma} \label{l4} The neighbors of a ($4,4)$-vertex are of degree at least 10.
\end{lemma}
\begin{proof} Suppose that there exists a $(4,4)$-vertex $v$ adjacent to a $9^-$-vertex. Let $G'$ be the graph obtained from $G$  after deleting $v$. $G'$ is proper with respect to $G$. By definition of minimal counterexample, we have $\chi_2(G') \leq 3\Delta +2.$ Since $d_2(v)  < 3\Delta+2$, we color $v$ by a safe color in $G$ to obtain $\chi_2(G) \leq 3\Delta +2$, a contradiction. \end{proof}
\begin{lemma} \label{l5} Let $v$ be a ($4$,$3$)-vertex. Then, we have the \mbox{following:}
\begin{enumerate}
    \item \label{l5p1} The neighbors of $v$ are of degree at least $8$.
    \item \label{l5p2} If $v$ is incident to a $4$-face, then the neighbors of $v$ are of degree at least $9$. 
\end{enumerate}
\end{lemma}
\begin{proof} 
Let $v$ be a ($4$,$3$)-vertex and set $N(v)=\{v_1,v_2,v_3,v_4\}$. W.l.o.g, suppose that the three 3-faces incident to $v$ are $vv_1v_2$, $vv_2v_3$ and $vv_3v_4$.
\begin{enumerate}
    \item Suppose that $v$ is adjacent to a $7^-$-vertex. Let $G'$ be the graph obtained from $G$ after deleting $v$ and adding edge $v_1v_4$. $G'$ is proper with respect to $G$. By definition of minimal counterexample, we have $\chi_2(G') \leq 3\Delta +2.$   Since $d_2(v)  < 3\Delta+2$, we color $v$ by a safe color in $G$ to obtain $\chi_2(G) \leq 3\Delta +2$, a contradiction.  
    \item Suppose that $v$ is incident to a $4$-face and suppose that $v$ is adjacent to an $8^-$-vertex.
Let $G'$ be the graph obtained from $G$  after deleting $v$. $G'$ is proper with respect to $G$. By definition of minimal counterexample, we have $\chi_2(G') \leq 3\Delta +2.$   Since $d_2(v)  < 3\Delta+2$, we color $v$ by a safe color in $G$ to obtain $\chi_2(G) \leq 3\Delta +2$, a contradiction. 
\end{enumerate}
\end{proof}
\begin{figure}[h]
\centering

\begin{minipage}{0.48\textwidth}
\centering
\begin{tikzpicture}[scale=1.2]

\coordinate (v1) at (-1,1);
\coordinate (v2) at (1,1);
\coordinate (v3) at (1,-1);
\coordinate (v4) at (-1,-1);
\coordinate (v)  at (0,0);

\draw (v1) -- (v2);
\draw (v1) -- (v);
\draw (v2) -- (v);
\draw (v3) -- (v);
\draw (v4) -- (v);
\draw (v3) -- (v4);

\fill (v1) circle (1pt);
\fill (v2) circle (1pt);
\fill (v3) circle (1pt);
\fill (v4) circle (1pt);
\fill (v)  circle (1pt);

\node[above] at (v1) {$v_1$};
\node[above] at (v2) {$v_2$};
\node[right] at (v3) {$v_3$};
\node[left]  at (v4) {$v_4$};
\node[right] at (v) {$v$};

\end{tikzpicture}
\end{minipage}
\hfill
\begin{minipage}{0.48\textwidth}
\centering
\begin{tikzpicture}[scale=1.2]

\coordinate (v1) at (-1,1);
\coordinate (v2) at (1,1);
\coordinate (v3) at (1,-1);
\coordinate (v4) at (-1,-1);
\coordinate (v)  at (0,0);

\draw (v1) -- (v2);
\draw (v1) -- (v);
\draw (v2) -- (v);
\draw (v3) -- (v);
\draw (v4) -- (v);
\draw (v2) -- (v3);

\fill (v1) circle (1pt);
\fill (v2) circle (1pt);
\fill (v3) circle (1pt);
\fill (v4) circle (1pt);
\fill (v)  circle (1pt);

\node[above] at (v1) {$v_1$};
\node[above] at (v2) {$v_2$};
\node[right] at (v3) {$v_3$};
\node[left]  at (v4) {$v_4$};
\node[right] at (v) {$v$};

\end{tikzpicture}
\end{minipage}

\caption{Illustration of Lemma 2.6}
\label{fig:two-graphs}
\end{figure}

\begin{lemma} \label{l6} Let $v$ be a ($4$,$2$)-vertex. Then, we have the following;
\begin{enumerate}
    \item \label{l6p1} The neighbors of $v$ are of degree at least $6$.
    \item \label{l6p2} If $v$ is adjacent to a $6$-vertex, then $v$ is a special vertex. 
\end{enumerate}
\end{lemma}
\begin{proof}
Let $v$ be a ($4$,$2$)-vertex and $f_1$ and $f_2$ be the two $3$-faces incident to $v$. W.l.o.g, suppose that $f_1=vv_1v_2$.
\begin{enumerate}
    \item Suppose that $v$ is adjacent to a $5^-$-vertex.
\begin{itemize}        
\item \textbf{Case 1:} $f_1$ and $f_2$ are adjacent: W.l.o.g, suppose that $f_2=vv_2v_3$. Let $G'$ be the graph obtained from $G$  after deleting $v$ and adding edge $v_2v_4$. $G'$ is proper with respect to $G$. By definition of minimal counterexample, we have $\chi_2(G') \leq 3\Delta +2.$   Since $d_2(v)  < 3\Delta+2$, we color $v$ by a safe color in $G$ to obtain $\chi_2(G) \leq 3\Delta +2$, a contradiction. 

\item \textbf{Case 2:} $f_1$ and $f_2$ are not adjacent:
Then, $f_2=vv_3v_4$. Let $G'$ be the graph obtained from $G$  after deleting $v$ and adding edges $v_1v_4$ and $v_2v_3$. $G'$ is proper with respect to $G$. By definition of minimal counterexample, we have $\chi_2(G') \leq 3\Delta +2.$   Since $d_2(v)  < 3\Delta+2$, we color $v$ by a safe color in $G$ to obtain $\chi_2(G) \leq 3\Delta +2$, a contradiction.
\end{itemize}
    \item Suppose that $v$ is adjacent to a $6$-vertex and suppose that there exists an edge in $G[N(v)]$ contained in two $3$-faces. 
  \begin{itemize}  
\item \textbf{Case 1:} $f_1$ and $f_2$ are adjacent: W.l.o.g, suppose that $f_2=vv_2v_3$. Let $G'$ be the graph obtained from $G$ after deleting $v$ and adding edge $v_2v_4$. $G'$ is proper with respect to $G$. By definition of minimal counterexample, we have $\chi_2(G') \leq 3\Delta +2.$   Since $d_2(v) < 3\Delta+2$, we color $v$ by a safe color in $G$ to obtain $\chi_2(G) \leq 3\Delta +2$, a contradiction. \\

\item \textbf{Case 2:} $f_1$ and $f_2$ are not adjacent:
Then, $f_2=vv_3v_4$. Let $G'$ be the graph obtained from $G$  after deleting $v$ and adding edges $v_1v_4$ and $v_2v_3$. $G'$ is proper with respect to $G$. By definition of minimal counterexample, we have $\chi_2(G') \leq 3\Delta +2.$   Since $d_2(v)  < 3\Delta+2$, we color $v$ by a safe color in $G$ to obtain $\chi_2(G) \leq 3\Delta +2$.
 \end{itemize}  
\end{enumerate}
\end{proof}

\begin{figure}[h]
\centering

\begin{minipage}{0.48\textwidth}
\centering
\begin{tikzpicture}[scale=1.2]

\coordinate (v1) at (-1,1);
\coordinate (v2) at (1,1);
\coordinate (v3) at (1,-1);
\coordinate (v4) at (-1,-1);
\coordinate (v)  at (0,0);

\coordinate (z) at (-2,0);
\coordinate (x) at (2,0);
\coordinate (y) at (0,-2);

\draw (v1) -- (v2);
\draw (v2) -- (x);
\draw (x) -- (v3);
\draw (v3) -- (y);
\draw (y) -- (v4);
\draw (v4) -- (z);
\draw (z) -- (v1);

\draw (v1) -- (v);
\draw (v2) -- (v);
\draw (v3) -- (v);
\draw (v4) -- (v);

\fill (v1) circle (1pt);
\fill (v2) circle (1pt);
\fill (v3) circle (1pt);
\fill (v4) circle (1pt);
\fill (v)  circle (1pt);
\fill (z) circle (1pt);
\fill (x) circle (1pt);
\fill (y) circle (1pt);

\node[above] at (v1) {$v_1$};
\node[above] at (v2) {$v_2$};
\node[right] at (v3) {$v_3$};
\node[left]  at (v4) {$v_4$};
\node[right] at (v) {$v$};
\node[left] at (z) {$z$};
\node[right] at (x) {$x$};
\node[below] at (y) {$y$};
\end{tikzpicture}
\end{minipage}
\hfill
\begin{minipage}{0.48\textwidth}
\centering
\begin{tikzpicture}[scale=1.2]

\coordinate (v1) at (-1,1);
\coordinate (v2) at (1,1);
\coordinate (v3) at (1,-1);
\coordinate (v4) at (-1,-1);
\coordinate (v)  at (0,0);

\coordinate (x) at (2,0);
\coordinate (y) at (0,-2);

\draw (v1) -- (v2);
\draw (v2) -- (x);
\draw (x) -- (v3);
\draw (v3) -- (y);
\draw (y) -- (v4);

\draw (v1) -- (v);
\draw (v2) -- (v);
\draw (v3) -- (v);
\draw (v4) -- (v);

\fill (v1) circle (1pt);
\fill (v2) circle (1pt);
\fill (v3) circle (1pt);
\fill (v4) circle (1pt);
\fill (v)  circle (1pt);
\fill (x) circle (1pt);
\fill (y) circle (1pt);

\node[above] at (v1) {$v_1$};
\node[above] at (v2) {$v_2$};
\node[right] at (v3) {$v_3$};
\node[left]  at (v4) {$v_4$};
\node[right] at (v) {$v$};
\node[right] at (x) {$x$};
\node[below] at (y) {$y$};
\end{tikzpicture}
\end{minipage}
\hfill
\begin{minipage}{0.48\textwidth}
\centering
\begin{tikzpicture}[scale=1.4]

\coordinate (v1) at (-1,1);
\coordinate (v2) at (1,1);
\coordinate (v3) at (1,-1);
\coordinate (v4) at (-1,-1);
\coordinate (v)  at (0,0);

\coordinate (y) at (-2,0);
\coordinate (x) at (2,0);

\draw (v1) -- (v2);
\draw (v2) -- (x);
\draw (x) -- (v3);

\draw (v4) -- (y);
\draw (y) -- (v1);

\draw (v1) -- (v);
\draw (v2) -- (v);
\draw (v3) -- (v);
\draw (v4) -- (v);

\fill (v1) circle (1pt);
\fill (v2) circle (1pt);
\fill (v3) circle (1pt);
\fill (v4) circle (1pt);
\fill (v)  circle (1pt);
\fill (y) circle (1pt);
\fill (x) circle (1pt);

\node[above] at (v1) {$v_1$};
\node[above] at (v2) {$v_2$};
\node[right] at (v3) {$v_3$};
\node[left]  at (v4) {$v_4$};
\node[right] at (v) {$v$};
\node[left] at (y) {$y$};
\node[right] at (x) {$x$};

\end{tikzpicture}
\end{minipage}
\caption{Illustration of Lemma 2.7}
\label{fig:two-graphs}
\end{figure}

\begin{lemma} \label{l7} Let $v$ be a ($4$,$1$)-vertex. Then, we have the following:
\begin{enumerate}
    \item \label{l7p1} If $v$ is incident to three $4$-faces, then $v$ is not adjacent to any $6^-$-vertex.
    \item \label{l7p2} If $v$ is incident to two $4$-faces, then $v$ is not adjacent to any $5^-$-vertex.
\end{enumerate}
\end{lemma}
\begin{proof}
Let $v$ be a ($4,1$)-vertex and $N(v)=\{v_1, v_2, v_3, v_4\}$. W.l.o.g, suppose that the $3$-face incident to $v$ is $vv_1v_2$.
\begin{enumerate}
\item Suppose that $v$ is a $(4,1,3)$-vertex  adjacent to a $6^-$-vertex.
W.l.o.g, suppose that the three $4$-faces incident  to $v$ are of the form $vv_2xv_3$,$vv_3yv_4$ and $vv_4zv_1$ where $x\in N(v_2)\cap N(v_3)$, $y\in N(v_3)\cap N(v_4)$ and $z\in N(v_1)\cap N(v_4)$.
Let $G’$ be the graph obtained from $G$ after deleting $v$ and adding edges $v_1v_4$ and $v_2v_3$. $G'$ is proper with respect to $G$. By definition of minimal counterexample, we have $\chi_2(G') \leq 3\Delta +2.$   Since $d_2(v)  < 3\Delta+2$, we color $v$ by a safe color in $G$ to obtain $\chi_2(G) \leq 3\Delta +2$, a contradiction.
\item  Let $v$ be a ($4,1,2$)-vertex and suppose to the contrary that $v$ is adjacent to a $5^-$-vertex. Let  $f_2$ and $f_3$ be the two $4$-faces incident to $v$. 
\begin{itemize}
\item \textbf{Case 1:} $f_2$ and $f_3$ are adjacent. W.l.o.g, suppose that $f_2=vv_2xv_3$ and $f_3=vv_3yv_4$ where $x\in N(v_2)\cap N(v_3)$ and $y\in N(v_3)\cap N(v_4)$. Let $G’$ be the graph obtained from $G$ after deleting $v$ and adding edges $v_2v_3$ and $v_1v_4$. $G'$ is proper with respect to $G$. By definition of minimal counterexample, we have $\chi_2(G') \leq 3\Delta +2.$   Since $d_2(v)  < 3\Delta+2$, we color $v$ by a safe color in $G$ to obtain $\chi_2(G) \leq 3\Delta +2$, a contradiction. \\

\item \textbf{Case 2:} $f_2$ and $f_3$ are not adjacent. W.l.o.g, suppose that $f_2=vv_2xv_3$ and $f_3=vv_4yv_1$ where $x\in N(v_2)\cap N(v_3)$ and $y\in N(v_1)\cap N(v_4)$.

\begin{itemize}
\item Subcase 1: $v_1$ or $v_2$ is a $5^-$ vertex: W.l.o.g, suppose that $v_1$ is a $5^-$-vertex.
Let $G’$ be the graph obtained from $G$ after deleting $v$ and adding edges $v_1v_3$ and $v_1v_4$. $G'$ is proper with respect to $G$. By definition of minimal counterexample, we have $\chi_2(G') \leq 3\Delta +2.$   Since $d_2(v)  < 3\Delta+2$, we color $v$ by a safe color in $G$ to obtain $\chi_2(G) \leq 3\Delta +2$, a contradiction. 
\item Subcase 2: $v_3$ or $v_4$ is a $5^-$ vertex: W.l.o.g, suppose that $v_4$ is a $5^-$-vertex. Let $G’$ be the graph obtained from $G$ after deleting $v$ and adding edges $v_1v_4$ and $v_3v_4$. $G'$ is proper with respect to $G$. By definition of minimal counterexample, we have $\chi_2(G') \leq 3\Delta +2.$   Since $d_2(v)  < 3\Delta+2$, \mbox{we color} $v$ by a safe color in $G$ to obtain $\chi_2(G) \leq 3\Delta +2$, a contradiction. 
\end{itemize}
\end{itemize}
\end{enumerate} 
\end{proof}
\begin{lemma} \label{l8} Let $v$ be a ($5,5$)-vertex. Then, we have the following:
\begin{enumerate}
\item \label{l8p1} If $v$ is adjacent to a $5^-$-vertex, then the other neighbors of $v$ are of degree at least 7.
\item \label{l8p2}  If $v$ is adjacent to a 5-vertex and a 7-vertex, then $v$ is a special vertex.
\item \label{l8p3} If $v$ is adjacent to two 6-vertices, then $v$ is a special vertex.
\end{enumerate}
\end{lemma}
\begin{proof} Let $v$ be a ($5,5$)-vertex.
\begin{enumerate}
\item Suppose that $v$ is adjacent to a $5^-$-vertex and another $6^-$-vertex. Let $G’$ be the graph obtained from $G$ after deleting $v$.
$G'$ is proper with respect to $G$. By definition of minimal counterexample, we have $\chi_2(G') \leq 3\Delta +2.$   Since $d_2(v)  < 3\Delta+2$, we color $v$ by a safe color in $G$ to obtain $\chi_2(G) \leq 3\Delta +2$, a contradiction. 
\item Suppose that $v$ is adjacent to a 5-vertex and a $7$-vertex and suppose to the contrary that there exists an edge in $G[N(v)]$ contained in two $3$-faces. Let $G’$ be the graph obtained from $G$ after deleting $v$. $G'$ is proper with respect to $G$. By definition of minimal counterexample, we have $\chi_2(G') \leq 3\Delta +2.$   Since $d_2(v)  < 3\Delta+2$, we color $v$ by a safe color in $G$ to obtain $\chi_2(G) \leq 3\Delta +2$, a contradiction. 
\item Suppose that $v$ is adjacent to two 6-vertices and suppose that there exists an edge in $G[N(v)]$ contained in two $3$-faces. Let $G’$ be the graph obtained from $G$ after deleting $v$.
$G'$ is proper with respect to $G$. By definition of minimal counterexample, we have $\chi_2(G') \leq 3\Delta +2.$   Since $d_2(v)  < 3\Delta+2$, we color $v$ by a safe color in $G$ to obtain $\chi_2(G) \leq 3\Delta +2$, a contradiction. 
\end{enumerate} 
\end{proof}
\begin{lemma} \label{l9} Let $v$ be a ($5,4,1$)-vertex.
Then, we have the following:
\begin{enumerate}
\item \label{l9p1} $v$ is adjacent to at most one $5^-$-vertex.
\item \label{l9p2} If $v$ is adjacent to a 6-vertex and a 5-vertex, then $v$ is a special vertex.
\item \label{l9p3} $v$ is adjacent to at most one ($6,6$)-vertex.
\end{enumerate}
\end{lemma}
\begin{proof} Let $v$ be a ($5,4,1$)-vertex and set $N(v)=\{v_1, v_2, v_3, v_4, v_5\}$. W.l.o.g, suppose that the four $3$-faces incident to $v$ are $vv_1v_2$, $vv_2v_3$, $vv_3v_4$ and $vv_4v_5$.
\begin{enumerate}
\item Suppose that $v$ is adjacent to two $5^-$-vertices. Let $G’$ be the graph obtained from $G$ after deleting $v$ and adding edge $v_1v_5$. $G'$ is proper with respect to $G$. By definition of minimal counterexample, we have $\chi_2(G') \leq 3\Delta +2.$   Since $d_2(v)  < 3\Delta+2$, we color $v$ by a safe color in $G$ to obtain $\chi_2(G) \leq 3\Delta +2$, a contradiction. 
\item Suppose that $v$ is adjacent to a 6-vertex and a 5-vertex and suppose that there exists an edge in $G[N(v)]$ contained in two $3$-faces.
Let $G’$ be the graph obtained from $G$ after deleting $v$ and adding edge $v_1v_5$. $G'$ is proper with respect to $G$. By definition of minimal counterexample, we have $\chi_2(G') \leq 3\Delta +2.$   Since $d_2(v)  < 3\Delta+2$, we color $v$ by a safe color in $G$ to obtain $\chi_2(G) \leq 3\Delta +2$, a contradiction. 
\item Suppose that $v$ has at least two ($6,6$)-neighbors. Then, at least three edges in $G[N(v)]$ are contained in two $3$-faces. Let $G’$ be the graph obtained from $G$ after deleting $v$ and adding edge $v_1v_5$. $G'$ is proper with respect to $G$. By definition of minimal counterexample, we have $\chi_2(G') \leq 3\Delta +2.$   Since $d_2(v)  < 3\Delta+2$, we color $v$ by a safe color in $G$ to obtain $\chi_2(G) \leq 3\Delta +2$, a contradiction.  
\end{enumerate}
 \end{proof}
\begin{lemma} \label{l10} Let $v$ be a ($5,4,0$)-vertex.
Then, we have the following:
\begin{enumerate}
\item \label{l10p1} If $v$ is adjacent to two $5^-$-vertices, then the other neighbors of $v$ are of degree at least $7$.
\item \label{l10p2} If $v$ is adjacent to two $5^-$-vertices and a $7$-vertex, then $v$ is a special vertex.
\item \label{l10p3} $v$ is adjacent to at most two ($6,6$)-vertices.
\item \label{l10p4} If $v$ is adjacent to a $5^-$-vertex, then $v$ is adjacent to at most one ($6,6$)-vertex.

\end{enumerate}
\end{lemma}
\begin{proof} Let $v$ be a ($5,4,0$)-vertex and set $N$($v$)=\{$v_1, v_2, v_3, v_4,v_5$\}. W.l.o.g, suppose that the four $3$-faces incident to $v$ are $vv_1v_2$, $vv_2v_3$, $vv_3v_4$ and $vv_4v_5$.
\begin{enumerate}
\item Suppose that $v$ is adjacent to two $5^-$-vertices and a $6$-vertex. Let $G'$ be the graph obtained from $G$ after deleting $v$ and adding edge $v_1v_5$. $G'$ is proper with respect to $G$. By definition of minimal counterexample, we have $\chi_2(G') \leq 3\Delta +2.$   Since $d_2(v)  < 3\Delta+2$, we color $v$ by a safe color in $G$ to obtain $\chi_2(G) \leq 3\Delta +2$, a contradiction. 
\item Suppose that $v$ is adjacent to two $5$-vertices and a $7$-vertex and suppose that there exists an edge in $G[N(v)]$ contained in two $3$-faces. Let $G’$ be the graph obtained from $G$ after deleting $v$ and adding edge $v_1v_5$. $G'$ is proper with respect to $G$. By definition of minimal counterexample, we have $\chi_2(G') \leq 3\Delta +2.$   Since $d_2(v)  < 3\Delta+2$, we color $v$ by a safe color in $G$ to obtain $\chi_2(G) \leq 3\Delta +2$, a contradiction. 
\item Suppose that $v$ is adjacent to three ($6,6$)-vertices. Then, at least four edges $G[N(v)]$ are contained in two $3$-faces. Let $G'$ be the graph obtained from $G$ after deleting $v$ and adding edge $v_1v_5$. $G'$ is proper with respect to $G$. By definition of minimal counterexample, we have $\chi_2(G') \leq 3\Delta +2.$   Since $d_2(v) \leq 2\Delta+6  < 3\Delta+2$, we color $v$ by a safe color in $G$ to obtain $\chi_2(G) \leq 3\Delta +2$, a contradiction.
\item Suppose that $v$ has a $5^-$-neighbor and suppose that $v$ is adjacent to two ($6,6$)-vertices. Then, at least three edges in $G[N(v)]$ are contained in two $3$-faces. Let $G'$ be the graph obtained from $G$  after deleting $v$ and adding edge $v_1v_5$. $G'$ is proper with respect to $G$. By definition of minimal counterexample, we have $\chi_2(G') \leq 3\Delta +2.$   Since $d_2(v)  < 3\Delta+2$, we color $v$ by a safe color in $G$ to obtain $\chi_2(G) \leq 3\Delta +2$, a contradiction. 
\end{enumerate}
\end{proof}
\begin{figure}[h]
\centering

\begin{tikzpicture}[scale=2, every node/.style={font=\small}]

\coordinate (v) at (0,0);

\coordinate (v1) at (120:1.1);
\coordinate (v2) at (60:1.1);
\coordinate (v3) at (0:1.1);
\coordinate (v4) at (-60:1.1);
\coordinate (v5) at (-120:1.1);
\coordinate (v6) at (180:1.1);
\coordinate (x) at (90:1.6);
\coordinate (y) at (30:1.7);
\coordinate (z) at (-30:1.6);
\coordinate (w) at (-90:1.7);

\fill (v1) circle (1pt);
\fill (v2) circle (1pt);
\fill (v3) circle (1pt);
\fill (v4) circle (1pt);
\fill (v5)  circle (1pt);
\fill (v6) circle (1pt);
\fill (v) circle (1pt);
\fill (x) circle (1pt);
\fill (y) circle (1pt);
\fill (z) circle (1pt);
\fill (w) circle (1pt);

\draw (v) -- (v2);
\draw (v) -- (v1);
\draw (v) -- (v3);
\draw (v) -- (v4);
\draw (v) -- (v5);
\draw (v) -- (v6);
\draw (v2) -- (v1);
\draw (v2) -- (v3);
\draw (v3) -- (v4);
\draw (v5) -- (v4);
\draw (v5) -- (v6);
\draw (x) -- (v1);
\draw (x) -- (v2);
\draw (y) -- (v2);
\draw (y) -- (v3);
\draw (z) -- (v4);
\draw (z) -- (v3);
\draw (w) -- (v4);
\draw (w) -- (v5);

\node[left]  at (v6) {$v_6$};
\node[above left]  at (v1) {$v_1$};
\node[above right] at (v2) {$v_2$};
\node[right] at (v3) {$v_3$};
\node[below right] at (v4) {$v_4$};
\node[below left] at (v5) {$v_5$};
\node[above right] at (v) {$v$};

\end{tikzpicture}
\caption{Illustration of Lemma 2.11}
\label{fig:two-graphs}
\end{figure}

\begin{lemma} \label{l11} Let $v$ be a ($6,5$)-vertex adjacent to two ($5,5$)-vertices. Then, $v$ is not adjacent to any ($5,4$)-vertex or any $4$-vertex. \end{lemma}
\begin{proof} Let $v$ be a ($6,5$)-vertex adjacent to two $(5,5)$-vertices. Then, at least four edges in $G[N(v)]$ are contained in two $3$-faces. Suppose to the contrary that $v$ is adjacent to a ($5,4$)-vertex or a $4$-vertex. Let $N$($v$)=\{$v_1, v_2, v_3, v_4, v_5, v_6$\}.  W.l.o.g, suppose that the five $3$-faces incident to $v$ are $vv_1v_2$, $vv_2v_3$, $vv_3v_4$, $vv_4v_5$ and $vv_5v_6$. Note that by \Cref{l8}(\ref{l8p1}), the ($5,5$)-vertices are not adjacent. Thus, either $v_2$ and $v_4$ are the $(5,5)$-vertices or $v_3$ and $v_5$ are the $(5,5)$-vertices. W.l.o.g, suppose that $v_2$ and $v_4$ are the ($5,5$)-vertices. Note that $v_6$ is the $(5,4)$-neighbor or the $4$-neighbor of $v$ since it is not adjacent to the $(5,5)$-vertices by \Cref{l8}(\ref{l8p1}). Note that if $v_6$ is a $(5,4)$-vertex, then five edges in $G[N(v)]$ are contained in two $3$-faces.
\begin{itemize}
\item \textbf{Case 1:} $\Delta =6$: Let $G'$ be the graph obtained from $G$ after deleting edge $vv_4$. $G'$ is proper with respect to $G$. By definition of minimal counterexample, we have $\chi_2(G') \leq 3\Delta +2=20$. Now, we uncolor $v$ and $v_4$ in $G$. Since $d_2(v) \leq 18$ and $d_2(v_4) \leq 18$, we color each vertex by a different safe color in $G$ and therefore we get $\chi_2(G) \leq 3\Delta +2$, a contradiction. 
\item \textbf{Case 2:} $\Delta \ge 7$: Let $G'$ be the graph obtained from $G$  after deleting $v$ and adding edges $v_1v_4$, $v_2v_4$, $v_4v_6$. $G'$ is proper with respect to $G$. By definition of minimal
counterexample, we have $\chi_2(G') \leq
3\Delta +2$. Since $d_2(v) \leq 2(\Delta -2) +(\Delta -1) +3+3+4-5 =3\Delta$, we color $v$ by a safe color in $G$ to obtain $\chi_2(G) \leq 3\Delta +2$, a contradiction. 
\end{itemize}
\end{proof}

\subsection{Charge Redistribution:}
Now, we will redistribute charges to prove that $G$ does not exist.\\
By Euler's formula, we have the following equality:
\\
$$\sum_{v\in V(G)}(d(v)-4)+ \sum_{f \in F(G) }(d(f)-4) = -8$$
Recall that an initial charge $d(v)-4$ is assigned to every vertex $v$ and $d(f)-4$ to every face $f$. \\
We redistribute charges among vertices and faces using the following discharging rules so that this redistribution yields nonnegative final charge to each vertex and each face:
\\
\textbf{R1:} Each $3$-face receives $\frac{1}{3}$ from every incident vertex.\\
\textbf{R2:} Each $5^+$-face sends $\frac{1}{3}$ to every incident $3$-vertex and $\frac{1}{5}$ to every incident $(\Delta-1)^-$-vertex. \\
\textbf{R3:} Each $3$-vertex receives $\frac{1}{9}$ from every neighbor.\\
\textbf{R4:} Each ($4,4$)-vertex receives $\frac{1}{3}$ from every neighbor. \\
\textbf{R5:} Each ($4,3$,$1$)-vertex receives $\frac{1}{4}$ from every neighbor. \\
\textbf{R6:} Each ($4,3,0$)-vertex receives $\frac{1}{5}$ from every neighbor.\\
\textbf{R7:} Each ($4,2$)-vertex receives $\frac{1}{6}$ from every neighbor.\\
\textbf{R8:} Each ($4,1,3$)-vertex receives $\frac{1}{12}$ from every neighbor. \\
\textbf{R9:} Each ($4,1,2$)-vertex receives $\frac{1}{30}$ from every neighbor. \\
\textbf{R10:} Each ($5$,$5$)-vertex receives $\frac{1}{6}$ from every $6^+$-neighbor except the $(6,6)$-vertex. \\
\textbf{R11:} Each ($5,4,1$)-vertex receives $\frac{1}{12}$ from every $6^+$-neighbor except the $(6,6)$-vertex.\\
\textbf{R12:} Each ($5,4,0$)-vertex receives $\frac{2}{45}$ from every $6^+$-neighbor except the $(6,6)$-vertex.\\

We call $v$ a bad $4$-vertex if $v$ is a $4$-vertex with negative charge after \mbox{applying R1 and R2.} We call $v$ a bad 5-vertex if $v$ is a 5-vertex with negative charge after \mbox{applying R1 and R2.} Note that the neighbors of a bad $4$-vertex are of degree at least 6 by \Cref{l4}, \ref{l5}, \ref{l6}, and \ref{l7}. Moreover, a $3$-vertex is not adjacent to any $5^-$-vertex by \Cref{l3}(\ref{l3p1}). We say $u$ is a bad $4$-neighbor of $v$ (resp., bad $5$-neighbor of $v$) if $u$ is a bad $4$-vertex adjacent to $v$ (resp., a bad $5$-vertex adjacent to $v$).
\\ 

Denote by $\mu(v) $ (resp., $\mu(f)$) the final charge of each vertex $v$ (resp., each face $f$). \\

Let $f\in F(G)$ and $v\in V(G)$. For each case of $f \in F(G)$ and $v\in V(G)$, we will prove that $\mu(v) \ge 0$ and $\mu(f) \ge 0$.
\begin{itemize}
\item If $f$ is a $3$-face: By R1, $f$ receives $\frac{1}{3}$ from each incident vertex. Then, $\mu (f)= d(f)-4+3.\frac{1}{3}= 0$.
\item If $f$ is a $4$-face: It does not send or receive any charge and so $\mu(f)=0$.
\item If $f$ is a $5$-face: $f$ is incident to at most two $3$-vertices
since $3$-vertices are not adjacent by \Cref{l3}(\ref{l3p1}). In this case, the remaining incident vertices are of degree $\Delta $ by \Cref{l3}(\ref{l3p1}). Hence, $f$ sends charge only to its incident $3$-vertices by R2. Then, $f$ sends charge most when all of its incident vertices are $(\Delta -1)^-$-vertices. Thus, we have $\mu (f) \ge d(f)-4-5.\frac{1}{5} \ge 0$.
\item If $f$ is a $6^+$-face: $f$ is incident to at most $\lfloor \frac{d(f)}{2}\rfloor$ $3$-vertices since $3$-vertices are not adjacent by \Cref{l3}(\ref{l3p1}). 
Note that $f$ sends charge the most when it has $\lfloor \frac{d(f)}{2}\rfloor$ incident $3$-vertices. Moreover, if $f$ has $\frac{d(f)}{2}$ incident $3$-vertices, then the remaining incident vertices to $f$ are $\Delta$-vertices. Hence, $f$ sends charge only to its incident $3$-vertices by R2. Then, $f$ sends charge most when all of its incident vertices are $(\Delta -1)^-$-vertices. Thus, we have $\mu (f) \ge d(f)-4- d(f).\frac{1}{5}  \ge 0$.
\item If $v$ is a $3$-vertex: \Cref{l3}(\ref{l3p3}), $v$ is incident to at least two $5^+$-faces. Hence, by $v$ receives $\frac{1}{3}$ from each incident $5^+$-face by R2 and $\frac{1}{3}$ from each neighbor by R3. Then, we have $\mu(v) \ge 3-4+2.\frac{1}{3}+ 3.\frac{1}{9} \ge 0$.
\item If $v$ is a $4$-vertex: Note that $v$ sends charge only to its incident $3$-faces if exists. Thus, if $v$ is not incident to any $3$-face, it doesn't send any charge and so we have $\mu (v) \ge d(v)-4\ge  0$. \\
If $v$ is a ($4$,$4$)-vertex, then $v$ sends $\frac{1}{3}$ to each incident $3$-face by R1 and receives $\frac{1}{3}$ from each neighbor by R4. Thus, $\mu (v)=4-4-4.\frac{1}{3}+4.\frac{1}{3}=0$. \\
If $v$ is a ($4$,$3$,1)-vertex, then $v$ sends $\frac{1}{3}$ to each incident $3$-face by R1 and receives $\frac{1}{4}$ from each neighbor by R5. Thus, $\mu (v)=4-4-3.\frac{1}{3}+4.\frac{1}{4}=0$. \\
If $v$ is a ($4$,$3$,$0$)-vertex, then $v$ sends $\frac{1}{3}$ to each incident $3$-face by R1 and receives $\frac{1}{5}$ from each neighbor by R6 and $\frac{1}{5}$ from its incident $5^+$-face by R3. Thus, we have $\mu (v)=4-4-3.\frac{1}{3}+4.\frac{1}{5}+\frac{1}{5}=0$. \\
If $v$ is a ($4$,$2$)-vertex, then $v$ sends $\frac{1}{3}$ to each incident $3$-face by R1 and receives $\frac{1}{6}$
from each neighbor by R7. Thus, $\mu (v)\ge 4-4-2.\frac{1}{3}+4.\frac{1}{6}\ge 0$. \\
If $v$ is a ($4$,$1$,$3$)-vertex, then $v$ sends $\frac{1}{3}$ to its incident 3-face by R1 and receives $\frac{1}{12}$
from each neighbor by R8. Thus,
$\mu (v)=4-4-\frac{1}{3}+4.\frac{1}{12}\ge0$. \\
If $v$ is a ($4$,$1$,$2$)-vertex, then $v$ sends $\frac{1}{3}$ to its incident $3$-face by R1 and receives $\frac{1}{30}$
from each neighbor by R9 and $\frac{1}{5}$ from its incident $5^+$-face.
Thus, we have $\mu (v)=4-4-\frac{1}{3}+4.\frac{1}{30}+\frac{1}{5} \ge0$. \\
If $v$ is a ($4$,$1$)-vertex incident to at most one $4$-face, then $v$ sends $\frac{1}{3}$ to its incident $3$-face and receives $\frac{1}{5}$
from at least two of its incident $5^+$-faces. Thus, we have $\mu (v)\ge 4-4-\frac{1}{3}+2.\frac{1}{5} \ge 0$.
\item If $v$ is a $5$-vertex: Note that $v$ is not adjacent to any bad $4$-vertex by \Cref{l4}, \ref{l5}, \ref{l6}(\ref{l6p1}) and \ref{l7}. Thus, $v$ sends charge only to its incident $3$-faces if exists. Hence, if $v$ is incident to at most three $3$-faces, we have $\mu (v)\ge 5-4-3.\frac{1}{3} \ge 0$. \\
If $v$ is a ($5$,$5$)-vertex: By Lemma 2.8(1), $v$ is adjacent to at most one $5^-$-vertex. If $v$ is adjacent to a $5^-$-vertex, then it is not adjacent to any other $6^-$-vertex by \Cref{l8}(\ref{l8p1}) and thus $v$ is not adjacent to any ($6$,$6$)-vertex. Then, $v$ receives $\frac{1}{6}$ from each $6^+$-neighbor  by R$10$ and sends $\frac{1}{3}$ to each incident $3$-face by R1. Thus, we have $\mu (v)=5-4-5.\frac{1}{3}+4.\frac{1}{6}=0$.
If $v$ is not adjacent to any $5^-$-vertex, we deduce that $v$ is adjacent to at most one ($6,6$)-vertex by \Cref{l8}(\ref{l8p3}). Thus, $v$ receives $\frac{1}{6}$ from at least four neighbors by R$10$ and sends $\frac{1}{3}$ to each incident $3$-face by R1. Then, we have $\mu (v)\ge 5-4-5.\frac{1}{3}+4.\frac{1}{6} \ge 0$. \\
If $v$ is a ($5,4,1$)-vertex: By \Cref{l9}(\ref{l9p1}), $v$ is adjacent to at most one $5^-$-vertex.
Suppose that $v$ is adjacent to a $5^-$-vertex. Then, we deduce that $v$ is not adjacent to any ($6,6$)-vertex by \Cref{l9}(\ref{l9p2}). Therefore, $v$ receives $\frac{1}{12}$ from each $6^+$-neighbor by R11 and sends $\frac{1}{3}$ to each incident $3$-face by R1. Hence, we have $\mu (v)=5-4-4.\frac{1}{3}+4.\frac{1}{12} =0$. Suppose now that $v$ is not adjacent to any $5^-$-vertex. Then, $v$ is adjacent to at most one ($6,6$)-vertex by \Cref{l9}(\ref{l9p3}). Thus, $v$ receives $\frac{1}{12}$ from at least four neighbors by R$11$ and sends $\frac{1}{3}$ to each incident $3$-face by R1. Then, we have $\mu (v)\ge 5-4-4.\frac{1}{3}+4.\frac{1}{12} \ge 0$. \\
If $v$ is a ($5,4,0$)-vertex: By R2, $v$ receives $\frac{1}{5}$ from its incident $5^+$-face since $\Delta \ge 6$. By \Cref{l10}(\ref{l10p1}), $v$ is adjacent to at most two $5^-$-vertices. Suppose that $v$ is adjacent to
two $5^-$-vertices. Then, $v$ is not adjacent to any ($6,6$)-vertex by \Cref{l10}(\ref{l10p1}). Hence, $v$ receives $\frac{2}{45}$ from three neighbors by R12 and sends $\frac{1}{3}$ to each incident $3$-face by R1. Therefore, we get $\mu (v)=5-4-4.\frac{1}{3}+\frac{1}{5}+3.\frac{2}{45}=0$.
Suppose that $v$ is adjacent to one $5^-$-vertex. Then, $v$ is adjacent to at most one ($6,6$)-vertex by \Cref{l10}(\ref{l10p4}). Hence, $v$ receives $\frac{2}{45}$ from at least three neighbors by R$12$ and sends $\frac{1}{3}$ to each incident $3$-face by R1. Thus, $\mu (v)\ge 5-4-4.\frac{1}{3}+\frac{1}{5}+3.\frac{2}{45} \ge 0$.
Suppose now that $v$ is not adjacent to any $5^-$-vertex. Then, $v$ is adjacent to at most two ($6,6$)-vertex by \Cref{l10}(\ref{l10p3}). Hence, $v$ receives $\frac{2}{45}$ from at least three neighbors by R$12$ and sends $\frac{1}{3}$ to each incident $3$-face by R1. Thus, $\mu (v)\ge 5-4-4.\frac{1}{3}+\frac{1}{5}+3.\frac{2}{45} \ge 0$.
\item If $v$ is a $6$-vertex: Note that $v$ is not adjacent to any ($4,1,3$)-vertex or ($4,3^+$)-vertex by \Cref{l4}, \ref{l5}(\ref{l5p1}), and \ref{l7}(\ref{l7p1}). Now, we will study the charge of $v$ according to number of $3$-faces incident to $v$. Recall that a bad $4$-vertex is not adjacent to any $5^-$-vertex.
\\
\textbf{1.} Suppose that $v$ is a ($6,6$)-vertex: Note that $v$ is not adjacent to any $(4,1)$-vertex since the neighbors of $v$ are incident to at least two $3$-faces. By \Cref{l4}, \ref{l5}(\ref{l5p1}), and \ref{l6}(\ref{l6p2}), $v$ is not adjacent to any $(4,2^+)$-vertex. By \Cref{l3}(\ref{l3p1}), $v$ is not adjacent to any $3$-vertex. Thus, $v$ does not send charge to any $4$-neighbor. Moreover, $v$ does not send charge to any bad $5$-neighbor by R10, R11 and R12.  Therefore, $v$ sends charge only to its incident $3$-faces. Hence, we have $\mu (v)= 6-4-6.\frac{1}{3} \ge 0$.
\\ \textbf{2.} Suppose that $v$ is ($6,5$)-vertex: Then, we deduce that the following properties about the neighbors of $v$:
\begin{itemize}
\item By \Cref{l3}(\ref{l3p2}), $v$ is not adjacent to any $3$-vertex.
\item By \Cref{l4}, \ref{l5}(\ref{l5p1}), \ref{l6}(\ref{l6p2}) and \ref{l7}(\ref{l7p1}), we deduce that the only bad $4$-vertices that could be adjacent to $v$ are the $(4,1,2)$-vertices. 
\item Since a $(4,1,2)$-vertex is incident to exactly one $3$-face, $v$ has at most two $(4,1,2)$-neighbors.
\item By \Cref{l8}(\ref{l8p1}), we deduce that if a $(5,5)$-vertex is adjacent to another $(5,5)$-vertex, it can not be adjacent to a $6$-vertex. Therefore, $v$ is adjacent to at most two ($5,5$)-vertices. 
\item If $v$ is adjacent to two $(5,5)$-vertices, it has no ($5,4^+$)-neighbor or $4$-neighbor by \Cref{l11}. Hence, in this case we have $\mu(v) = 6-4-5.\frac{1}{3}-2.\frac{1}{6} \ge 0.$ 
\item If $v$ is adjacent to one ($5,5$)-vertex, we deduce that it is adjacent to at most two ($5,4$)-vertices by \Cref{l8}(\ref{l8p1}) and \ref{l10}(\ref{l10p1}). Hence,  we have $\mu(v) \ge 6-4-5.\frac{1}{3}-\frac{1}{6} -2.\frac{1}{12}\ge 0.$
\item If $v$ is not adjacent to any ($5,5$)-vertex, we deduce that it is adjacent to at most four ($5,4$)-vertices. Hence, we have $\mu(v) \ge 6-4-5.\frac{1}{3}-4.\frac{1}{12} \ge 0.$
\end{itemize}
Note that $v$ sends charge to its bad 5-neighbor more than it sends to its bad $4$-neighbor and so the worst case occurs when $v$ has bad $5$-neighbors. Thus, we get $\mu (v) \ge 0$. \\
\textbf{3.} Suppose that $v$ is a ($6,4^-$)-vertex: 
Since $v$ sends charge to its incident $3$-face more than it sends to any neighbor, the worst case occurs when $v$ is a ($6,4$)-vertex. Hence, we may suppose that $v$ is a ($6,4$)-vertex. We have following properties about the neighbors of $v$:
\begin{itemize}
    \item The only bad $4$-vertices that could be adjacent to $v$ are the ($4,2$)-vertices and the ($4,1,2$)-vertices by \Cref{l4}, \ref{l5}, and \ref{l7}(\ref{l7p1}).
     Thus, $v$ sends charge to its bad 5-neighbor more than it sends to its bad $4$-neighbor except the $(4,2)$-neighbor. Hence, the worst case occurs when $v$ has no $(4,1,2)$-neighbors. 
    \item Since bad $4$-vertices are not adjacent, $v$ has at most four $(4,2)$-neighbors. Since the neighbors of a $(4,2)$-vertex are of degree at least $6$ by \Cref{l6}(\ref{l6p1}), we deduce that in this case $v$ has no other bad neighbor. Thus, $v$ sends $\frac{1}{3}$ to each incident $3$-face by R1 and $\frac{1}{6}$ to each bad neighbor by R7. Hence, $\mu(v)=6-4-4.\frac{1}{3}-4.\frac{1}{6}\ge 0.$
    \item Since a $3$-vertex is not incident to any $3$-face by \Cref{l3}(\ref{l3p2}), $v$ has at most one $3$-neighbor.
    \item By \Cref{l8}(\ref{l8p1}), we deduce that if a $(5,5)$-vertex is adjacent to another $(5,5)$-vertex, it can not be adjacent to a $6$-vertex. Therefore, $v$ is adjacent to at most two ($5,5$)-vertices.
    \item If $v$ has two ($5,5$)-neighbors, it has no other ($5,4$)-neighbor since a $(5,5)$ can not be adjacent to both a $6$-vertex and a $5$-vertex by \Cref{l8}(\ref{l8p1}). Then, $v$ sends $\frac{1}{3}$ to each incident $3$-face by R1, $\frac{1}{6}$ to each ($5,5$)-neighbor by R12 and at most $\frac{1}{6}$ to its sixth bad neighbor if exists. Thus, we have $\mu(v) \ge 6-4-4.\frac{1}{3}-2.\frac{1}{6}-\frac{1}{6}\ge 0.$
    \item If $v$ has one ($5,5$)-neighbor, it has at most two ($4,2$)-neighbors which is the worst case since the neighbors of a bad $4$-vertex are of degree at least 6. Then, $v$ sends $\frac{1}{3}$ to each incident $3$-face by R1, $\frac{1}{6}$ to its $(5,5)$-neighbor by R10 and $\frac{1}{6}$ to each $(4,2)$-neighbor by R7. Thus, we have $\mu(v) \ge 6-4-4.\frac{1}{3}-\frac{1}{6}-2.\frac{1}{6}\ge 0.$
    \item If $v$ has no ($5,5$)-neighbor, we deduce that it has at most four ($4,2$)-neighbors since bad $4$-vertices are not adjacent. Then, $v$ sends $\frac{1}{3}$ to each incident $3$-face by R1 and $\frac{1}{6}$ to each $(4,2)$-neighbor by R7. Thus, we have $\mu(v) \ge 6-4-4.\frac{1}{3}-4.\frac{1}{6}\ge 0.$
\end{itemize}
In all cases, we get that $\mu(v)  \ge 0$.
\item If $v$ is a $7$-vertex: The worst case occurs when $v$ is a ($7$,$7$)-vertex since $v$ sends charge to its incident $3$-faces more than it sends to any other neighbor. Hence, we suppose that $v$ is a ($7,7$)-vertex. Then, since the neighbors of $v$ are incident to at least two $3$-faces, we deduce that $v$ is not adjacent to any $(4,1)$-vertex. By \Cref{l4} and \ref{l5}(\ref{l5p1}), we deduce that the only bad $4$-vertices that could be adjacent to $v$ are the $(4,2)$-vertices. Therefore, since $(4,2)$-vertices are not adjacent, the worst case occurs when the bad neighbors of $v$ are the ($5,5$)-vertices and the $(5,4,1)$-vertices. By \Cref{l10}(\ref{l10p2}), we deduce that $v$ has at most four ($5,4$)-neighbors and in this case $v$ is not adjacent to any other bad $4$-vertex or bad $5$-vertex. If a ($5,5$)-vertex is adjacent to another ($5,5$)-vertex, it can't be adjacent to a $(7,7)$-vertex by \Cref{l8}(\ref{l8p2}). Therefore, $v$ has at most three ($5,5$)-neighbors. Thus, the worst case occurs when $v$ has three ($5,5$)-neighbors and in this case $v$ is not adjacent to any other bad $4$-vertex or bad $5$-vertex by \Cref{l8}(\ref{l8p2}). Hence, $v$ sends $\frac{1}{3}$ to each incident $3$-face by R1 and $\frac{1}{6}$ to each bad neighbor by R10. Thus, we have $\mu (v) = 7-4-7.\frac{1}{3} -3.\frac{1}{6} \ge 0$.
\item If $v$ is an $8$-vertex: By \Cref{l4} and \ref{l5}(\ref{l5p1}) Lemma 2.4 and 2.5(1), $v$ is not adjacent to any $(4,4)$-vertex or $(4,3,1)$-vertex. By \Cref{l8}(\ref{l8p1}), we deduce that $v$ has at most five $(5,5)$-neighbors and in this case $v$ is not adjacent to any other bad $4$-vertex or bad $5$-vertex. Note that the ($4,3$)-vertices are not adjacent by \Cref{l5}(\ref{l5p1}). Thus, the worst case occurs when $v$ is an ($8,8$)-vertex and has four $(4,3,0)$-neighbors. Hence, $v$ sends $\frac{1}{3}$ to each incident $3$-face by R1 and $\frac{1}{5}$ to each bad neighbor by R6. Thus, we get $\mu (v)= 8-4-8.\frac{1}{3}-4.\frac{1}{5} \ge 0$.
\item If $v$ is a $9$-vertex: By \Cref{l4}, $v$ is not adjacent to any $(4,4)$-vertex. By \Cref{l8}(\ref{l8p1}), we deduce that $v$ has at most six $(5,5)$-neighbors and in this case $v$ is not adjacent to any other bad $4$-vertex or bad $5$-vertex. Note that the ($4,3,1$)-vertices are not adjacent by \Cref{l5}(\ref{l5p2}). Thus, the worst case occurs when $v$ is an ($9$,$9$)-vertex and has four ($4$,$3$,$1$)-neighbors. Hence, $v$ sends $\frac{1}{3}$ to each incident $3$-face by R1 and $\frac{1}{4}$ to each bad neighbor by R5. Thus, we get $\mu (v)\ge 9-4-9.\frac{1}{3}-4.\frac{1}{4} \ge 0$.
\item If $v$ is a $k$-vertex such that $k \ge 10$: By \Cref{l8}(\ref{l8p1}), we deduce that $v$ has at most $\lfloor\frac{2k}{3}\rfloor$ $(5,5)$-neighbors and in this case $v$ is not adjacent to any other bad $4$-vertex or bad $5$-vertex. Note that ($4,4$)-vertices are not adjacent by \Cref{l4}. Moreover, $v$ sends charge the most to its adjacent $(4,4)$-vertices if exists. Thus, the worst case occurs when $v$ is a ($k$,$k$)-vertex having $\lfloor\frac{k}{2} \rfloor$ $(4,4)$-neighbors and in such a case $v$ has no other bad neighbor by \Cref{l4}. Hence, in this case $v$ sends $\frac{1}{3}$ to each incident $3$-face by R1 and $\frac{1}{3}$ to each bad neighbor by R4.  Thus, we get $\mu (v) \ge k-4-\frac{k}{3}-\frac{k}{2}.\frac{1}{3} \ge 0$ for $k \ge 10$.

\end{itemize}

\vspace{0.5cm}
\textbf{Acknowledgment:} I would like to express my sincere gratitude to Dr. Maydoun Mortada for her support and guidance in the process of writing this paper.

\end{large}

\end{document}